\documentclass[11pt]{amsart}
\usepackage[margin=1.2in]{geometry}

\usepackage{doi}
\usepackage{amsmath, amstext, amsbsy, amsopn, amsfonts,amssymb}
\usepackage{array,multirow}
\usepackage{tikz}
\usetikzlibrary{graphs}
\usetikzlibrary{graphs.standard}

\DeclareMathOperator\conv {conv} 
\DeclareMathOperator\TSP {TSP} 
\DeclareMathOperator\ATSP {ATSP}

\usepackage{comment}
\usepackage{algorithm}
\usepackage{algpseudocode}

\theoremstyle{plain}
\newtheorem{theorem}{Theorem}
\newtheorem{lemma}{Lemma}

\begin{document}

\title[Hamiltonian decomposition and the 1-skeleton of the TSP polytope by VNS]{Hamiltonian decomposition and verifying vertex adjacency in 1-skeleton of the traveling salesperson polytope by variable neighborhood search}

\author{Andrei Nikolaev \and Anna Kozlova}
\thanks {The research is supported by the grant of the President of the Russian Federation MK-2620.2018.1 (agreement no. 075-015-2019-746).}

\address{%
	P.\,G. Demidov Yaroslavl State University, Yaroslavl, Russia
}
\email{andrei.v.nikolaev@gmail.com, fyz95@mail.ru}

\begin{abstract}
	We consider a Hamiltonian decomposition problem of partitioning a regular graph into edge-disjoint Hamiltonian cycles. 
	A sufficient condition for vertex adjacency in the 1-skeleton of the traveling salesperson polytope can be formulated as the Hamiltonian decomposition problem in a 4-regular multigraph.
	We introduce a heuristic general variable neighborhood search algorithm for this problem based on finding a vertex-disjoint cycle cover of the multigraph through reduction to perfect matching and several cycle merging operations.
	The algorithm has a one-sided error: the answer ``not adjacent'' is always correct, and was tested on random directed and undirected Hamiltonian cycles and on pyramidal tours.
\end{abstract}

\maketitle

\section{Introduction}

\textit{A Hamiltonian decomposition} of a regular graph is a partition of its edge set into Hamiltonian cycles. The problem of finding edge-disjoint Hamiltonian cycles in a given regular graph plays an important role in combinatorial optimization \cite{krar:1995}, coding theory \cite{bae:bose:2003,bail:2009}, privacy-preserving distributed mining algorithms \cite{clif:kant:vaid:2002}, analysis of interconnection networks \cite{hung:2011} and other areas. 
See also theoretical results on estimating the number of Hamiltonian decompositions of regular graphs \cite{gleb:lur:sud:2017}. 
Our motivation for this problem comes from the field of polyhedral combinatorics.

\section{Traveling salesperson polytope}

We consider a classic traveling salesperson problem: given a complete weighted graph (or digraph) $K_n=(V,E)$, it is required to find a Hamiltonian cycle of minimum weight.
We denote by $HC_{n}$ the set of all Hamiltonian cycles in $K_{n}$.
With each Hamiltonian cycle $x \in HC_n$ we associate a characteristic vector $x^v \in \mathbb{R}^{E}$ by the following rule:
\[
x^v_e = 
\begin{cases}
1,& \text{ if the cycle } x \text{ contains an edge } e \in E,\\
0,& \text{ otherwise. }
\end{cases}
\]
The polytope
\[\TSP(n) = \conv \{x^v \ | \ x \in HC_n \}\]
is called \textit{the symmetric traveling salesperson polytope}.

\textit{The asymmetric traveling salesperson polytope} $\ATSP(n)$ is defined similarly as the convex hull of characteristic vectors of all possible Hamiltonian cycles in the complete digraph $K_{n}$.

The 1-\textit{skeleton} of a polytope $P$ is the graph whose vertex set is the vertex set of $P$ and edge set is the set of geometric edges or one-dimensional faces of $P$. 
The study of 1-skeleton is of interest, since, on the one hand, there are algorithms for perfect matching, set covering, independent set, a ranking of objects, problems with fuzzy measures, and many others that are based on the vertex adjacency relation in 1-skeleton and the local search technique (see, for example, \cite{ag:katz:tol:2017,bali:1985,cher:ham:1987,comb:mir:2010,gab:1977}). On the other hand, some characteristics of 1-skeleton, such as the diameter and clique number, estimate the time complexity for different computation models and classes of algorithms \cite{bond:1983,bond:nik:2013,grot:padb:1985}.

Unfortunately, the classic result by Papadimitriou states that the construction of 1-skeleton of the traveling salesperson polytope is hard for both directed and undirected graphs.

\begin{theorem} [Papadimitriou \cite{papa:1978}]\label{theorem:papa}
	The question of whether two vertices of the polytopes $\TSP(n)$ or $\ATSP(n)$ are nonadjacent is NP-complete. 
\end{theorem}

As a result, there are a lot of papers on the diameter and clique number of 1-skeleton of $\TSP(n)$ and $\ATSP(n)$ \cite{bond:1983,risp:1998,sier:1993}, but little progress with adjacency relation.
We can only note the polynomial-time algorithms to test vertex adjacencies in the pedigree polytope \cite{arth:2006}, the polytope of pyramidal tours \cite{bond:nik:2018} and the polytope of pyramidal tours with step-backs \cite{nik:2019} which are directly related to the traveling salesperson problem.

\section{Hamiltonian decomposition and the 1-skeleton}

We apply the Hamiltonian decomposition problem to analyze the 1-skeleton of the traveling salesperson polytope.
Let $x = (V,E_x)$ and $y=(V,E_y)$ be two Hamiltonian cycles on the vertex set $V$.
We denote by $x \cup y$ a multigraph $(V,E_x \cup E_y)$ that contains all edges of both cycles $x$ and $y$.

\begin{lemma}[Sufficient condition for nonadjacency \cite{rao:1976}]\label{lemma_sufficient}
	Given two Hamiltonian cycles $x$ and $y$, if the multigraph $x \cup y$ contains a Hamiltonian decomposition into edge-disjoint cycles $z$ and $w$ different from $x$ and $y$, then the corresponding vertices $x^v$ and $y^v$ of the polytope $\TSP(n)$ (or $\ATSP(n)$) are not adjacent.
\end{lemma}

\begin{figure}[t]
	\centering
	\begin{tikzpicture}[scale=0.85]
	\begin{scope}[every node/.style={circle,thick,draw}]
	\node (1) at (0,0) {1};
	\node (2) at (1,1) {2};
	\node (3) at (2.5,1) {3};
	\node (4) at (3.5,0) {4};
	\node (5) at (2.5,-1) {5};
	\node (6) at (1,-1) {6};
	\end{scope}
	\draw [line width=0.3mm] (1) edge (2);
	\draw [line width=0.3mm] (2) edge (3);
	\draw [line width=0.3mm] (3) edge (4);
	\draw [line width=0.3mm] (4) edge (5);
	\draw [line width=0.3mm] (5) edge (6);
	\draw [line width=0.3mm] (6) edge (1);
	\draw (1.75, -1.7) node{\textit{x}};
	
	\begin{scope}[yshift=-4cm]
	\begin{scope}[every node/.style={circle,thick,draw}]
	\node (1) at (0,0) {1};
	\node (2) at (1,1) {2};
	\node (3) at (2.5,1) {3};
	\node (4) at (3.5,0) {4};
	\node (5) at (2.5,-1) {5};
	\node (6) at (1,-1) {6};
	\end{scope}
	\draw [line width=0.3mm] (1) edge (4);
	\draw [line width=0.3mm] (4) edge (6);
	\draw [line width=0.3mm] (6) edge (2);
	\draw [line width=0.3mm] (2) edge (3);
	\draw [line width=0.3mm] (3) edge (5);
	\draw [line width=0.3mm] (5) edge (1);
	\draw (1.75, -1.7) node{\textit{y}};
	\end{scope}
	
	\begin{scope}[xshift=4.5cm,yshift=-2cm]
	\begin{scope}[every node/.style={circle,thick,draw}]
	\node (1) at (0,0) {1};
	\node (2) at (1,1) {2};
	\node (3) at (2.5,1) {3};
	\node (4) at (3.5,0) {4};
	\node (5) at (2.5,-1) {5};
	\node (6) at (1,-1) {6};
	\end{scope}
	\draw [line width=0.3mm] (1) edge (2);
	\draw [line width=0.3mm, bend right=10] (2) edge (3);
	\draw [line width=0.3mm] (3) edge (4);
	\draw [line width=0.3mm] (4) edge (5);
	\draw [line width=0.3mm] (5) edge (6);
	\draw [line width=0.3mm] (6) edge (1);
	\draw [line width=0.3mm] (1) edge (4);
	\draw [line width=0.3mm] (4) edge (6);
	\draw [line width=0.3mm] (6) edge (2);
	\draw [line width=0.3mm, bend left=10] (2) edge (3);
	\draw [line width=0.3mm] (3) edge (5);
	\draw [line width=0.3mm] (5) edge (1);	
	\draw (1.75, -1.7) node{\textit{x $\cup$ y}};
	\end{scope}

	\begin{scope}[xshift=9cm]
	\begin{scope}[every node/.style={circle,thick,draw}]
	\node (1) at (0,0) {1};
	\node (2) at (1,1) {2};
	\node (3) at (2.5,1) {3};
	\node (4) at (3.5,0) {4};
	\node (5) at (2.5,-1) {5};
	\node (6) at (1,-1) {6};
	\end{scope}
	\draw [line width=0.3mm] (1) edge (4);
	\draw [line width=0.3mm] (4) edge (5);
	\draw [line width=0.3mm] (5) edge (3);
	\draw [line width=0.3mm] (3) edge (2);
	\draw [line width=0.3mm] (2) edge (6);
	\draw [line width=0.3mm] (6) edge (1);
	\draw (1.75, -1.7) node{\textit{z}};
	
	\begin{scope}[yshift=-4cm]
	\begin{scope}[every node/.style={circle,thick,draw}]
	\node (1) at (0,0) {1};
	\node (2) at (1,1) {2};
	\node (3) at (2.5,1) {3};
	\node (4) at (3.5,0) {4};
	\node (5) at (2.5,-1) {5};
	\node (6) at (1,-1) {6};
	\end{scope}
	\draw [line width=0.3mm] (1) edge (2);
	\draw [line width=0.3mm] (2) edge (3);
	\draw [line width=0.3mm] (3) edge (4);
	\draw [line width=0.3mm] (4) edge (6);
	\draw [line width=0.3mm] (6) edge (5);
	\draw [line width=0.3mm] (5) edge (1);
	\draw (1.75, -1.7) node{\textit{w}};
	\end{scope}
	\end{scope}
	\end{tikzpicture}
	\caption{The multigraph $x \cup y$ has two different Hamiltonian decompositions}
	\label{image:not_adjacent}
\end{figure}
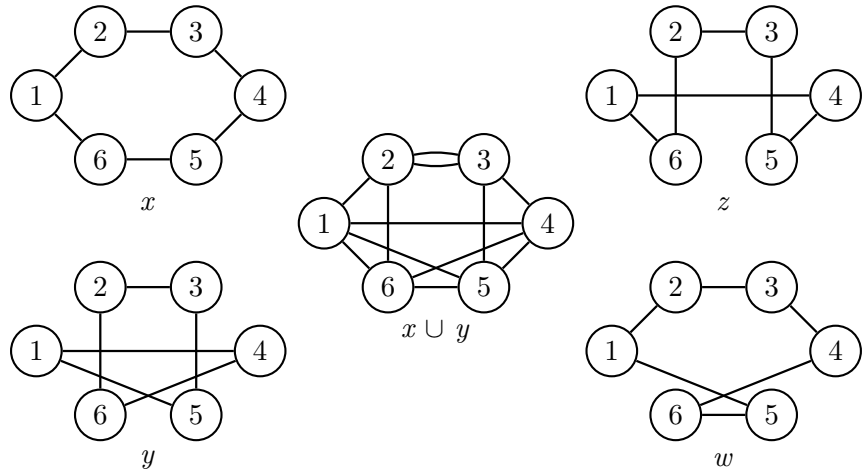

From a geometric point of view, the sufficient condition means that the segment connecting two vertices $x^v$ and $y^v$ intersects with the segment connecting two other vertices $z^v$ and $w^v$ of the polytope $\TSP(n)$ (or $\ATSP(n)$ correspondingly). Thus, the vertices $x^v$ and $y^v$ cannot be adjacent in 1-skeleton. An example of a satisfied sufficient condition is shown in Fig.~\ref{image:not_adjacent}.

Thus, the sufficient condition for vertex nonadjacency of the traveling salesperson polytope can be formulated as a combinatorial problem.

\vspace{1mm}
\textsc{Instance.} Let $x$ and $y$ be two Hamiltonian cycles.

\textsc{Question.} Does the multigraph $x \cup y$ contain a pair of edge-disjoint Hamiltonian cycles $z$ and $w$ different from $x$ and $y$ such that
\[z \cup w = x \cup y \text{ and } z \cap w = \emptyset?\]

Testing of whether a graph has a Hamiltonian decomposition is NP-complete, even for 4-regular undirected graphs and 2-regular directed graphs \cite{per:1984}. Thus, verifying the sufficient condition for vertex nonadjacency is hard.
Polynomially solvable special cases of the Hamiltonian decomposition and vertex adjacency problems have been studied in the literature, in particular: pedigrees \cite{arth:2006}, pyramidal tours \cite{bond:nik:2018} and pyramidal tours with step-backs \cite{nik:2019}.
Based on them, we can consider polynomial sufficient conditions for nonadjacency in the traveling salesperson polytope, however, all of them will be weaker than the condition of Lemma~\ref{lemma_sufficient} \cite{arth:2013}.

Note that any 4-regular multigraph always has an even number of Hamiltonian decompositions \cite{thom:1978}. 
There is no contradiction with Lemma~\ref{lemma_sufficient} and Theorem~\ref{theorem:papa}.
If the Hamiltonian cycles $x$ and $y$ contain the same edge, then both copies of this edge are included in the multigraph $x \cup y$.
The exchange of such edges can provide a Hamiltonian decomposition, however, the cycles $z$ and $w$ will not differ from $x$ and $y$.

\section{General variable neighborhood search}

We introduce a heuristic general variable neighborhood search algorithm with 3 different neighborhood structures to find the Hamiltonian decomposition in the multigraph $x \cup y$.
The variable neighborhood search metaheuristic was proposed by Mladenović and Hansen in 1997 \cite{mlad:1997} and has since developed rapidly both in its methods and applications. See, for example, current surveys \cite{han:mlad:2018,han:mlad:tod:han:2017} and some recent applications of general VNS \cite{law:luc:mlad:2019,mlad:del:lap:wil:2020}.

The GVNS algorithm is an improved version of the simulated annealing approach from \cite{kozl:nik:2019} that has the same feasible set, the objective function, and the third neighborhood structure.

\subsection{Feasible set}

\begin{figure}[t]
	\centering
	\begin{tikzpicture}[scale=0.85]
	\begin{scope}[every node/.style={circle,thick,draw}]
	\node (a1) at (0,0) {1};
	\node (a2) at (1,1) {2};
	\node (a3) at (2.5,1) {3};
	\node (a4) at (3.5,0) {4};
	\node (a5) at (2.5,-1) {5};
	\node (a6) at (1,-1) {6};
	\end{scope}
	
	\node at (1.75, -1.75) {$x \cup y$};
	
	\draw [thick] (a1) edge (a2);
	\draw [thick] (a1) edge (a6);
	\draw [thick] (a2) edge (a6);
	\draw [thick] (a3) edge (a4);
	\draw [thick] (a3) edge (a5);
	\draw [thick] (a4) edge (a5);
	
	\draw [thick, bend left=25] (a1) edge (a2);
	\draw [thick] (a1) edge (a3);
	\draw [thick] (a2) edge (a3);
	\draw [thick] (a4) edge (a6);
	\draw [thick] (a5) edge (a6);
	\draw [thick, bend left=25] (a4) edge (a5);

	\begin{scope}[xshift=-4.75cm]
	\begin{scope}[every node/.style={circle,thick,draw}]
	\node (b1) at (0,0) {1};
	\node (b2) at (1,1) {2};
	\node (b3) at (2.5,1) {3};
	\node (b4) at (3.5,0) {4};
	\node (b5) at (2.5,-1) {5};
	\node (b6) at (1,-1) {6};
	\end{scope}
	
	\draw [thick] (b1) edge (b2);
	\draw [thick] (b1) edge (b6);
	\draw [thick] (b2) edge (b6);
	\draw [thick] (b3) edge (b4);
	\draw [thick] (b3) edge (b5);
	\draw [thick] (b4) edge (b5);
	\node at (1.75, -1.75) {$z$};

	\end{scope}

	\begin{scope}[xshift=4.75cm]
	\begin{scope}[every node/.style={circle,thick,draw}]
	\node (c1) at (0,0) {1};
	\node (c2) at (1,1) {2};
	\node (c3) at (2.5,1) {3};
	\node (c4) at (3.5,0) {4};
	\node (c5) at (2.5,-1) {5};
	\node (c6) at (1,-1) {6};
	\end{scope}
	
	\node at (1.75, -1.75) {$w = (x \cup y) \backslash z$};
	
	\draw [thick, bend left=25] (c1) edge (c2);
	\draw [thick] (c1) edge (c3);
	\draw [thick] (c2) edge (c3);
	\draw [thick] (c4) edge (c6);
	\draw [thick] (c5) edge (c6);
	\draw [thick, bend left=25] (c4) edge (c5);	
	\end{scope}
	\end{tikzpicture}
	\caption{The multigraph $x \cup y$ and its two complementary cycle covers}
	\label{image:cycle_covers}
\end{figure}
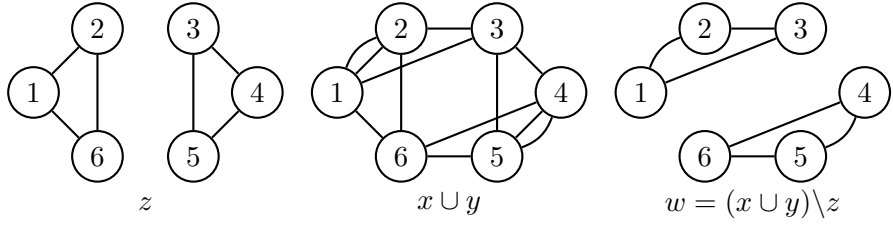

The set of feasible solutions consists of pairs $z$ and $w$ of vertex-disjoint cycle covers of the multigraph $x \cup y$.
\textit{A vertex-disjoint cycle cover} of a graph $G$ is a set of cycles with no vertices in common which are subgraphs of $G$ and contain all vertices of $G$.
The idea is as follows. All vertices in the multigraph $x \cup y$ have degrees equal to $4$ (or both indegrees and outdegrees equal to $2$ for directed cycles). Let $z$ be a vertex-disjoint cycle cover of $x \cup y$, then all the remaining edges form a graph $w = (x \cup y) \backslash z$ with all vertex degrees being equal to $2$ (both indegrees and outdegrees equal to $1$). Thus, $w$ is also a vertex-disjoint cycle cover of $x \cup y$ (Fig.~\ref{image:cycle_covers}).

Finding a vertex-disjoint cycle cover of both directed and undirected graphs can be performed in polynomial time by a reduction to perfect matching \cite{tutt:1954}.
Let us recall that \textit{a perfect matching} is a set of pairwise nonadjacent edges that matches all vertices of the graph.
The procedures for directed and undirected graphs are somewhat different. We consider them separately. 

\begin{figure}[t]
	\centering
	\begin{tikzpicture}[scale=0.8]
	\begin{scope}[every node/.style={circle,thick,draw}]
	\node (A) at (0,0) {$v$};
	\end{scope}
	\draw [thick] (A) -- node [right] {$A$} (0,2);
	\draw [thick] (A) -- node [below] {$B$} (2,0);
	\draw [thick] (A) -- node [left] {$C$} (0,-2);
	\draw [thick] (A) -- node [above] {$D$} (-2,0);
	\draw (0, -2.5) node{vertex $v$ of $G$};
	\node at (2.75,0) {$\Rightarrow$};
	\begin{scope}[xshift=7cm]
	\begin{scope}[every node/.style={circle,thick,draw}]
	\node (A) at (0,2) {$v_a$};
	\node (B) at (2,0) {$v_b$};
	\node (C) at (0,-2) {$v_c$};
	\node (D) at (-2,0) {$v_d$};
	\node (E) at (0.5,0.5) {$v_1$};
	\node (F) at (-0.5,-0.5) {$v_2$};
	\end{scope}
	\draw [thick] (E) edge (A);
	\draw [thick] (E) edge (B);
	\draw [thick] (E) edge (C);
	\draw [thick] (E) edge (D);
	\draw [thick] (F) edge (A);
	\draw [thick] (F) edge (B);
	\draw [thick] (F) edge (C);
	\draw [thick] (F) edge (D);
	\draw [thick] (A) -- node [right] {$A$} (0,3.5);
	\draw [thick] (B) -- node [below] {$B$} (3.5,0);
	\draw [thick] (C) -- node [left] {$C$} (0,-3.5);
	\draw [thick] (D) -- node [above] {$D$} (-3.5,0);
	\draw (0, -4) node{gadget $G_v$ of $G'$};
	\end{scope}
	\end{tikzpicture}
	\caption{Construction of the graph $G'$ for the undirected multigraph $G$}
	\label{image:gadget_undirected}
\end{figure}
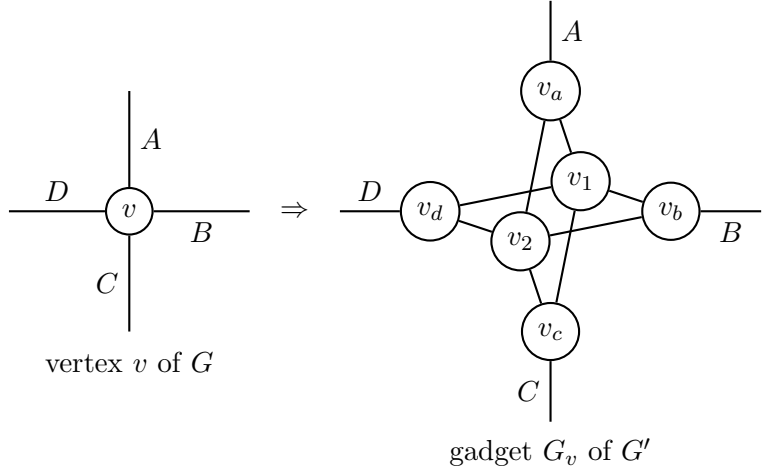

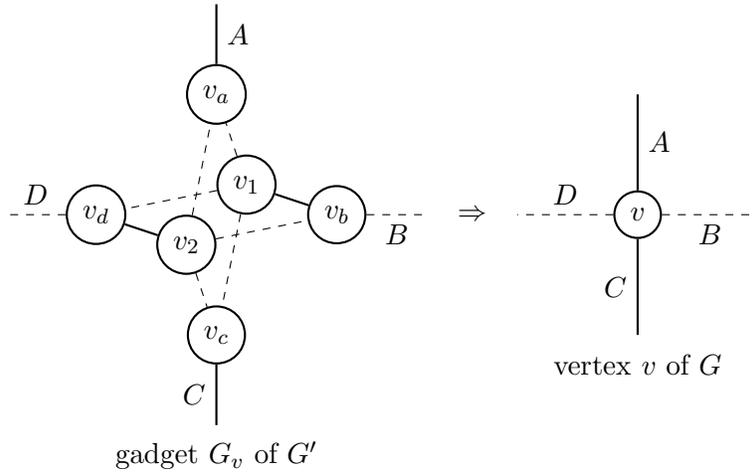
\begin{figure}[t]
	\centering
	\begin{tikzpicture}[scale=0.8]
	\begin{scope}[every node/.style={circle,thick,draw}]
	\node (A) at (0,2) {$v_a$};
	\node (B) at (2,0) {$v_b$};
	\node (C) at (0,-2) {$v_c$};
	\node (D) at (-2,0) {$v_d$};
	\node (E) at (0.5,0.5) {$v_1$};
	\node (F) at (-0.5,-0.5) {$v_2$};
	\end{scope}
	\draw [dashed] (E) edge (A);
	\draw [thick] (E) edge (B);
	\draw [dashed] (E) edge (C);
	\draw [dashed] (E) edge (D);
	\draw [dashed] (F) edge (A);
	\draw [dashed] (F) edge (B);
	\draw [dashed] (F) edge (C);
	\draw [thick] (F) edge (D);
	\draw [thick] (A) -- node [right] {$A$} (0,3.5);
	\draw [dashed] (B) -- node [below] {$B$} (3.5,0);
	\draw [thick] (C) -- node [left] {$C$} (0,-3.5);
	\draw [dashed] (D) -- node [above] {$D$} (-3.5,0);
	\draw (0, -4) node{gadget $G_v$ of $G'$};
	
	\node at (4.25,0) {$\Rightarrow$};
	
	\begin{scope}[xshift=7cm]
	\begin{scope}[every node/.style={circle,thick,draw}]
	\node (A) at (0,0) {$v$};
	\end{scope}
	\draw [thick] (A) -- node [right] {$A$} (0,2);
	\draw [dashed] (A) -- node [below] {$B$} (2,0);
	\draw [thick] (A) -- node [left] {$C$} (0,-2);
	\draw [dashed] (A) -- node [above] {$D$} (-2,0);
	\draw (0, -2.5) node{vertex $v$ of $G$};
	
	\end{scope}
	\end{tikzpicture}
	\caption{Perfect matching in $G'$ (solid) and vertex-disjoint cycle cover in $G$}
	\label{image:matching_undirected}
\end{figure}

Let $x$ and $y$ be two undirected Hamiltonian cycles.
\begin{enumerate}
	\item [Step 1.] From the multigraph $x \cup y = G = (V, E)$, we construct a new graph $G' = (V', E')$. With each vertex $v \in V$ we associate a gadget $G_v$ that is a complete bipartite subgraph $K_{4,2}$ (note that the degree of $v$ equals 4) as it is shown in Fig.~\ref{image:gadget_undirected}:
	\begin{itemize}
		\item there are $4$ vertices in the outer part ($v_a$, $v_b$, $v_c$ and $v_d$) that correspond to $4$ edges incident to $v$ in $G$ (edges $A$, $B$, $C$, $D$) and are connected with other gadgets;
		\item there are $2$ vertices in the inner part ($v_1$ and $v_2$) that are only adjacent to the vertices inside the gadget.
	\end{itemize}
	\item [Step 2.] A perfect matching in $G'$ corresponds to a vertex-disjoint cycle cover in the original graph $G$.
	Indeed, perfect matching has to cover both vertices in the inner part. 
	Therefore, it includes exactly one edge adjacent to $v_1$ and one edge adjacent to $v_2$.
	These two edges cover only two of the four vertices in the outer part $\{v_a,v_b,v_c,v_d\}$.
	Two other vertices have to be covered by the edges of $G$ (Fig.~\ref{image:matching_undirected}). 
	We include these edges into $z$, then the degree of each vertex $v$ in the graph $z$ equals $2$, and therefore, $z$ is a vertex-disjoint cycle cover of the multigraph $x \cup y$.
\end{enumerate}

A perfect matching in a general undirected graph can be found by Micali-Vazirani algorithm in $O(\sqrt{V} E)$ time \cite{mic:1980}.

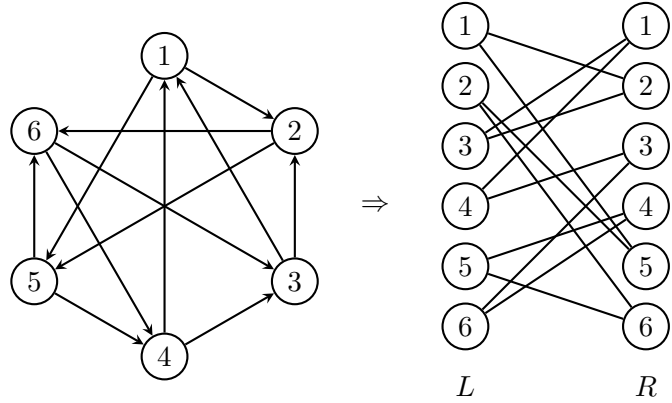
\begin{figure}[t]
	\centering
	\begin{tikzpicture}[scale=0.8]
	\begin{scope}[every node/.style={circle,thick,draw,inner sep=3pt}]
	\graph [clockwise] {
		subgraph I_n [n=6,name=A, radius=2cm]; 
		
		A 1 --[->,>=stealth,thick] A 2;
		A 1 --[->,>=stealth,thick] A 5; 
		A 2 --[->,>=stealth,thick] A 5;
		A 2 --[->,>=stealth,thick] A 6;
		A 3 --[->,>=stealth,thick] A 1;
		A 3 --[->,>=stealth,thick] A 2;
		A 4 --[->,>=stealth,thick] A 1;
		A 4 --[->,>=stealth,thick] A 3;
		A 5 --[->,>=stealth,thick] A 4;
		A 5 --[->,>=stealth,thick] A 6;
		A 6 --[->,>=stealth,thick] A 3;
		A 6 --[->,>=stealth,thick] A 4;
	};
	\end{scope}
	
	\node at (3.5,0) {$\Rightarrow$};
	
	\begin{scope}[xshift=5cm,yshift=-2cm]
	\begin{scope}[every node/.style={circle,thick,draw,inner sep=3pt}]
	\node (A1) at (0,5) {1};
	\node (A2) at (0,4) {2};
	\node (A3) at (0,3) {3};
	\node (A4) at (0,2) {4};
	\node (A5) at (0,1) {5};
	\node (A6) at (0,0) {6};
	\node (B1) at (3,5) {1};
	\node (B2) at (3,4) {2};
	\node (B3) at (3,3) {3};
	\node (B4) at (3,2) {4};
	\node (B5) at (3,1) {5};
	\node (B6) at (3,0) {6};
	\end{scope}
	
	\node at (0,-1) {$L$};
	\node at (3,-1) {$R$};
	
	\draw [thick] (A1) edge (B2);
	\draw [thick] (A1) edge (B5);
	\draw [thick] (A2) edge (B5);
	\draw [thick] (A2) edge (B6);
	\draw [thick] (A3) edge (B1);
	\draw [thick] (A3) edge (B2);
	\draw [thick] (A4) edge (B1);
	\draw [thick] (A4) edge (B3);
	\draw [thick] (A5) edge (B4);
	\draw [thick] (A5) edge (B6);
	\draw [thick] (A6) edge (B3);
	\draw [thick] (A6) edge (B4);
	\end{scope}
	\end{tikzpicture}
	\caption{Directed multigraph $D$ (left) and corresponding bipartite graph $D'$ (right)}
	\label{image:gadget_directed}
\end{figure}

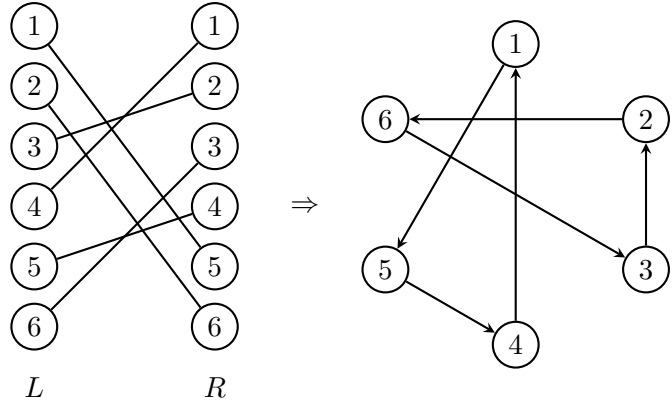
\begin{figure}[t]
	\centering
	\begin{tikzpicture}[scale=0.8]
	\begin{scope}[every node/.style={circle,thick,draw,inner sep=3pt}]
	\node (A1) at (0,5) {1};
	\node (A2) at (0,4) {2};
	\node (A3) at (0,3) {3};
	\node (A4) at (0,2) {4};
	\node (A5) at (0,1) {5};
	\node (A6) at (0,0) {6};
	\node (B1) at (3,5) {1};
	\node (B2) at (3,4) {2};
	\node (B3) at (3,3) {3};
	\node (B4) at (3,2) {4};
	\node (B5) at (3,1) {5};
	\node (B6) at (3,0) {6};
	\end{scope}
	
	\node at (0,-1) {$L$};
	\node at (3,-1) {$R$};
	
	\draw [thick] (A1) edge (B5);
	\draw [thick] (A2) edge (B6);
	\draw [thick] (A3) edge (B2);
	\draw [thick] (A4) edge (B1);
	\draw [thick] (A5) edge (B4);
	\draw [thick] (A6) edge (B3);
	
	\node at (4.5,2) {$\Rightarrow$};
	
	\begin{scope}[xshift=8cm,yshift=2.2cm]
	\begin{scope}[every node/.style={circle,thick,draw,inner sep=3pt}]
	\graph [clockwise] {
		subgraph I_n [n=6,name=A, radius=2cm]; 
		
		A 1 --[->,>=stealth,thick] A 5; 
		A 2 --[->,>=stealth,thick] A 6;
		A 3 --[->,>=stealth,thick] A 2;
		A 4 --[->,>=stealth,thick] A 1;
		A 5 --[->,>=stealth,thick] A 4;
		A 6 --[->,>=stealth,thick] A 3;
	};
	\end{scope}
	\end{scope}
	\end{tikzpicture}
	\caption{Perfect matching in $D'$ (left) and vertex-disjoint cycle cover of $D$ (right)}
	\label{image:matching_directed}
\end{figure}

Let $x$ and $y$ be two directed Hamiltonian tours.
\begin{enumerate}
	\item [Step 1.] From the directed multigraph $x \cup y = D = (V, A)$, we construct a bipartite graph $D' = (L, R, E)$. With each vertex $v \in V$ we associate a pair of vertices $v_L \in L$ and $v_R \in R$, and with each edge $(u,v) \in A$ we associate a new edge $(u_L, v_R)$ in the bipartite graph $D'$ (Fig.~\ref{image:gadget_directed}).
	\item [Step 2.] A perfect matching in the bipartite graph $D'$ corresponds to a vertex-disjoint directed cycle cover in the original graph $D$.
	Indeed, every vertex of $D$ is a head of exactly one edge and a tail of exactly one edge of a perfect matching in $D'$ (Fig.~\ref{image:matching_directed}).
\end{enumerate}

A perfect matching in a bipartite graph can be found by Hopcroft–Karp algorithm in $O(\sqrt{V} E)$ time \cite{hop:1973}.

\subsection{Objective function}

As the objective function, we chose the total number of connected components in the vertex-disjoint cycle covers $z$ and $w$.
If it equals 2, then $z$ and $w$ are Hamiltonian cycles.

\subsection{First neighborhood structure}

We consider the directed graphs.
Without loss of generality, let here and further $z$ be a vertex-disjoint cycle cover with at least 2 connected components.

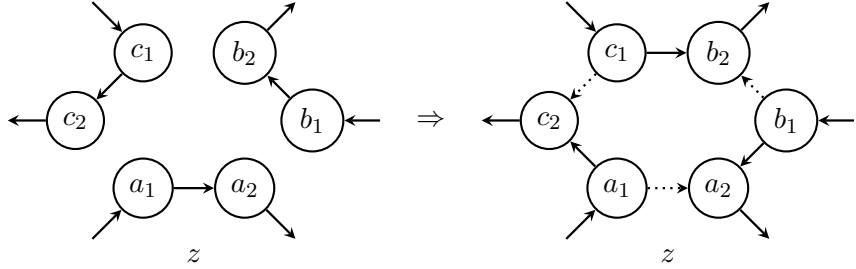
\begin{figure}[t]
	\centering
	\begin{tikzpicture}[scale=0.9]
	\begin{scope}[every node/.style={circle,thick,draw}]
	\node (A1) at (0,0) {$a_1$};
	\node (A2) at (1.5,0) {$a_2$};
	\node (B1) at (2.5,1) {$b_1$};
	\node (B2) at (1.5,2) {$b_2$};
	\node (C1) at (0,2) {$c_1$};
	\node (C2) at (-1,1) {$c_2$};
	\end{scope}
	\draw [->,>=stealth,line width=0.3mm] (-0.75,-0.75) -- (A1);
	\draw [->,>=stealth,line width=0.3mm] (A1) edge (A2);
	\draw [->,>=stealth,line width=0.3mm] (A2) -- (2.25,-0.75);
	\draw [->,>=stealth,line width=0.3mm] (B2) -- (2.25,2.75);
	\draw [->,>=stealth,line width=0.3mm] (B1) edge (B2);
	\draw [->,>=stealth,line width=0.3mm] (3.5,1) -- (B1);
	\draw [->,>=stealth,line width=0.3mm] (-0.75,2.75) -- (C1);
	\draw [->,>=stealth,line width=0.3mm] (C1) edge (C2);
	\draw [->,>=stealth,line width=0.3mm] (C2) -- (-2,1);
	
	\node at (0.75,-1) {$z$};
	
	\node at (4.25,1) {$\Rightarrow$};
	
	\begin{scope}[xshift=7cm]
	\begin{scope}[every node/.style={circle,thick,draw}]
	\node (A1) at (0,0) {$a_1$};
	\node (A2) at (1.5,0) {$a_2$};
	\node (B1) at (2.5,1) {$b_1$};
	\node (B2) at (1.5,2) {$b_2$};
	\node (C1) at (0,2) {$c_1$};
	\node (C2) at (-1,1) {$c_2$};
	\end{scope}
	\draw [->,>=stealth,line width=0.3mm] (-0.75,-0.75) -- (A1);
	\draw [->,>=stealth,line width=0.3mm,dotted] (A1) edge (A2);
	\draw [->,>=stealth,line width=0.3mm] (A2) -- (2.25,-0.75);
	\draw [->,>=stealth,line width=0.3mm] (B2) -- (2.25,2.75);
	\draw [->,>=stealth,line width=0.3mm,dotted] (B1) edge (B2);
	\draw [->,>=stealth,line width=0.3mm] (B1) edge (A2);
	\draw [->,>=stealth,line width=0.3mm] (3.5,1) -- (B1);
	\draw [->,>=stealth,line width=0.3mm] (-0.75,2.75) -- (C1);
	\draw [->,>=stealth,line width=0.3mm,dotted] (C1) edge (C2);
	\draw [->,>=stealth,line width=0.3mm] (C1) edge (B2);
	\draw [->,>=stealth,line width=0.3mm] (C2) -- (-2,1);
	\draw [->,>=stealth,line width=0.3mm] (A1) edge (C2);
	
	\node at (0.75,-1) {$z$};
	\end{scope}
	\end{tikzpicture}
	\caption{Example of the first neighborhood structure}
	\label{Fig_first_neighborhood_structure}
\end{figure}

The idea is as follows. We choose an edge $(a_1,a_2)$ of $z$ and remove it (move to $w$).
Since the directed multigraph $x \cup y$ is 2-regular, there are only 2 edges with their heads adjacent to the vertex $a_2$: $(a_1,a_2)$ and $(b_1,a_2)$. 
We add the edge $(b_1,a_2)$ of $w$ to $z$ (Fig.~\ref{Fig_first_neighborhood_structure}, the removed edges are dotted). 
But now outdegree of $b_1$ in $z$ equals 2. Again, there are only 2 edges with their tails adjacent to $b_1$: $(b_1,a_2)$ and $(b_1,b_2)$. We cannot backtrack and remove the just added edge $(b_1,a_2)$. 
Therefore, we remove the edge $(b_1,b_2)$, and so on. This is summarized in Algorithm \ref{Alg:LS_1}.

\begin{algorithm}[t]
	\caption{Local search with respect to the first neighborhood}\label{Alg:LS_1}
	\begin{algorithmic}[1]
		\Procedure{LocalSearch\_1}{$z,w$}	
		\While {there are unchecked and unfixed edges in $z$ (or $w$)}
		\State $z' \gets z$, $w' \gets w$ 		\Comment{Neighboring solution}
		\State Choose an unchecked and unfixed edge $(a_1,a_2)$ of $z'$ 	\Comment{In random order}
		\State $z' \gg (a_1,a_2) \gg w'$ 		\Comment{Move the edge $(a_1,a_2)$ from $z'$ to $w'$}
		\State $z' \ll (b_1,a_2) \ll w'$ 		\Comment{Move the edge $(b_1,a_2)$ from $w'$ to $z'$}
		\State $a_1.outdegree \gets a_1.outdegree - 1$ 
		\State $b_1.outdegree \gets b_1.outdegree + 1$ 
		\While {there is some vertex $b_1$ of $z'$ with outdegree equals 2}  	\Comment{See Fig.~\ref{Fig_first_neighborhood_structure}}
		\State $z' \gg (b_1,b_2) \gg w'$
		\State $b_1.outdegree \gets b_1.outdegree - 1$ 
		\State $z' \ll (c_1,b_2) \ll w'$
		\State $c_1.outdegree \gets c_1.outdegree + 1$ 
		\EndWhile
		\If {a better solution is found}
		\State \Return $z'$ and $w'$	\Comment{First improvement w.r.t. the first neighborhood}
		\EndIf
		\State Mark all moved edges as checked	\Comment{Avoid double checking the same edges}
		\EndWhile
		\State \Return $z$ and $w$ is a local minimum w.r.t. the first neighborhood
		\EndProcedure
	\end{algorithmic}
\end{algorithm}

Note that the whole procedure is deterministic.
If after this operation $z$ and $w$ are the correct cycle covers and the number of connected components has decreased, then we proceed to a new solution.
The total size of the neighborhood is $O(V)$, as we have to check every edge of the vertex-disjoint cycle cover $z$ exactly once.

We call an edge \textit{fixed} if it coincides in the original cycles $x$ and $y$. Since two copies of one edge cannot get into one cycle, we fix them in $z$ and $w$.

For the undirected graphs, the procedure is similar to Algorithm \ref{Alg:LS_1}. 
Again in a cycle for each edge of $z$ we move the edge to $w$ and try to restore the cycle covers by exchanging edges.
The key difference is that the procedure is no longer deterministic: at each step, we have a choice of several edges to add or remove (see the second neighborhood structure for more details). Thus, every time we make a random choice of which edge to add or remove and run $k$ random walks.
The time complexity of checking the first neighborhood for undirected graphs is $O(k \cdot V^2)$.

\subsection{Second neighborhood structure}

Now we consider the undirected graphs. Let again $z$ be the cycle cover with at least two connected components. We choose an edge of $z$ and remove it (move to $w$). 
Then in $z$ there are vertices with degrees other than 2.

\begin{figure}[t]
	\centering
	\begin{tikzpicture}[scale=0.95]
	
	\begin{scope}[every node/.style={circle,draw,thick,inner sep=3pt}]
	\node (N) at (0,1.25) {$v_1$};
	\node (C) at (0,0) {$~u~$};
	\node (S) at (0,-1.25) {$v_3$};
	\node (W) at (-1.25,0) {$v_4$};
	\node (E) at (1.25,0) {$v_2$};
	\end{scope}
	
	\draw [thick,dotted] (N) -- (C);
	\draw [thick] (C) -- (S);

	\draw [->,>=stealth,thick] (-1.25,-1.8) -- (-2.5,-2.5);
	\draw [->,>=stealth,thick] (0,-1.8) -- (0,-2.5);
	\draw [->,>=stealth,thick] (1.25,-1.8) -- (2.5,-2.5);
	
	\begin{scope}[xshift=-3.75cm,yshift=-4.24cm]
	\begin{scope}[every node/.style={circle,draw,thick,inner sep=3pt}]
	\node (N) at (0,1.25) {$v_1$};
	\node (C) at (0,0) {$~u~$};
	\node (S) at (0,-1.25) {$v_3$};
	\node (W) at (-1.25,0) {$v_4$};
	\node (E) at (1.25,0) {$v_2$};
	\end{scope}
	
	\draw [thick,dotted] (N) -- (C);
	\draw [thick] (C) -- (S);
	\draw [thick] (C) -- (E);
	\end{scope}
	
	\begin{scope}[yshift=-4.25cm]
	\begin{scope}[every node/.style={circle,draw,thick,inner sep=3pt}]
	\node (N) at (0,1.25) {$v_1$};
	\node (C) at (0,0) {$~u~$};
	\node (S) at (0,-1.25) {$v_3$};
	\node (W) at (-1.25,0) {$v_4$};
	\node (E) at (1.25,0) {$v_2$};
	\end{scope}
	
	\draw [thick,dotted] (N) -- (C);
	\draw [thick] (C) -- (S);
	\draw [thick] (C) -- (W);
	\end{scope}

	\begin{scope}[xshift=3.75cm,yshift=-4.25cm]
	\begin{scope}[every node/.style={circle,draw,thick,inner sep=3pt}]
	\node (N) at (0,1.25) {$v_1$};
	\node (C) at (0,0) {$~u~$};
	\node (S) at (0,-1.25) {$v_3$};
	\node (W) at (-1.25,0) {$v_4$};
	\node (E) at (1.25,0) {$v_2$};
	\end{scope}
	
	\draw [thick,dotted] (N) -- (C);
	\draw [thick,dotted] (C) -- (S);
	\draw [thick] (C) -- (W);
	\draw [thick] (C) -- (E);
	\end{scope}
	
	\end{tikzpicture}
	\caption {Second neighborhood structure, 3 ways to restore the vertex degree}
	\label{Fig_second_neighborhood_structure}
\end{figure}
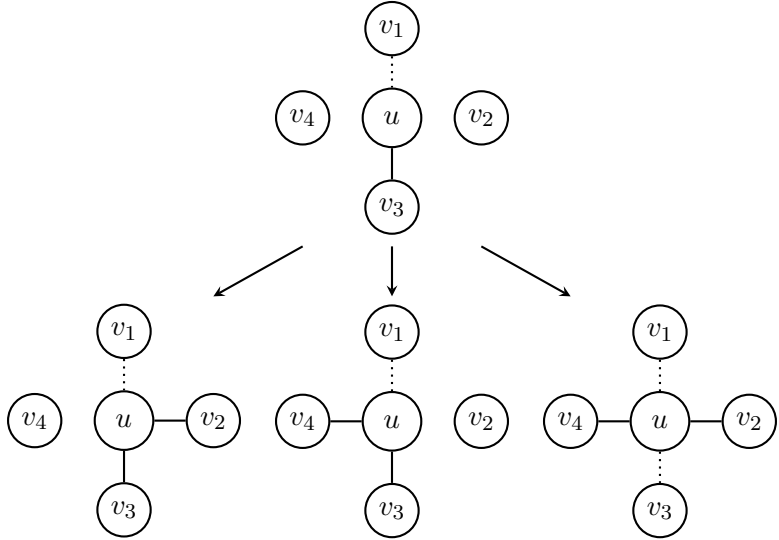

Let $u$ be some vertex of degree $1$.
We can restore the vertex degree of $u$ in three different ways: by adding one of two incident edges from $w$, or by adding both incident edges from $w$ at the same time and removing the existing edge of $z$ (Fig.~\ref{Fig_second_neighborhood_structure}, the removed edges are dotted).
Note that we cannot return the previously deleted edge $(u,v_1)$.
We explore all three options using DFS with a bounded search tree.
Again, if we manage to merge cycles, get the correct cycle covers, and reduce the number of connected components, then we proceed to a new solution.
The procedure is summarized in Algorithm \ref{Alg:second_neighborhood_structure}.

\begin{algorithm}[t]
	\caption{Local search with respect to the second neighborhood}\label{Alg:second_neighborhood_structure}
	\begin{algorithmic}[1]
		\Procedure{LocalSearch\_2}{$z,w,depthL$}
		\While {there are unchecked and unfixed edges in $z$ (or $w$)}
		\State Choose an unchecked and unfixed edge $(u,v_1)$ of $z$	\Comment{In random order}
		\State \Call{BoundedSearchTree}{$z \gg (u,v_1) \gg w,depth=1,depthL$} \Comment{Move the edge $(u,v_1)$ from $z$ to $w$ and explore the search tree on modified graphs}
		\If {a better solution is found}
		\State \Return $z$ and $w$ 		\Comment{First improvement w.r.t. the second neighborhood}
		\EndIf
		\State Mark the edge $(u,v_1)$ as checked
		\EndWhile
		\State \Return $z$ and $w$ is a local minimum w.r.t. the second neighborhood
		\EndProcedure
		\vspace{2mm}
		\Procedure{BoundedSearchTree}{$z \text{ and } w,depth,depthL$}
		\If {$z$ and $w$ are correct cycle covers}
		\If {the number of connected components has decreased}
		\State Proceed to a new solution and exit recursion
		\EndIf
		\State \Return				\Comment{Backtrack to the previous step}
		\EndIf
		\If {$depth > depthL$ or other stop conditions are satisfied}
		\State \Return				\Comment{Backtrack to the previous step}
		\EndIf
		\State Choose a vertex $u$ of $z$ whose degree is not equal to 2
		\If {vertex degree of $u$ is equal to 1}  	\Comment{See Fig.~\ref{Fig_second_neighborhood_structure}}
		\State \Call{BoundedSearchTree}{$z \ll (u,v_2) \ll w,depth+1,depthL$}  
		\State \Call{BoundedSearchTree}{$z \ll (u,v_4) \ll w, depth+1,depthL$}
		\State \Call{BoundedSearchTree}{$z \ll (u,v_2),(u,v_4) \ll w \text{ and } z \gg (u,v_3) \gg w, depth + 1,depthL$}
		\EndIf 
		\State Consider other cases of vertex degrees
		\EndProcedure
		
	\end{algorithmic}
\end{algorithm}

Other cases of vertex degree are treated similarly. If a vertex degree is equal to 3, we again have up to three ways to restore a vertex degree by removing one or two edges. In cases of degrees 0 and 4, there is only one way to restore a vertex degree by adding/removing two edges at the same time.

The time complexity of checking the second neighborhood is $O(3^{d} \cdot V)$ where $d = depthL$ is the recursion depth limit.
The procedure for directed graphs is similar, but with no more than two options at each step.

Note that the first and second neighborhood structures explore the same search tree for merging cycles through the exchange of edges.
However, the difference is as follows:
\begin{itemize}
	\item the first neighborhood structure explores only some branches (branches of a special type for directed graphs and a limited number of random branches for undirected graphs);
	\item the second neighborhood structure explores the search tree completely, but only to a limited depth.
\end{itemize}

\subsection{Third neighborhood structure}

In the third neighborhood structure, we choose a random edge of $w$ with endpoints in two different connected components of $z$ (or vice versa) and add it to the queue of fixed edges.
Such an edge always exists due to the connectivity of the multigraph $x \cup y$.
The neighboring solution is constructed as a new vertex-disjoint cycle cover by the perfect matching algorithms with fixed edges forming the initial matching.
When the queue limit $fixEdgesN$ is exceeded, the first edge of the queue is deleted.
The procedure is summarized in Algorithm \ref{Alg:SA}.

\begin{algorithm}[t]
	\caption{Simulated annealing with respect to the third neighborhood}\label{Alg:SA}
	\begin{algorithmic}[1]
		\Procedure{SimulatedAnnealing}{$z,w,iter,T,fixEdgesN$}
		\State Choose a random edge $e$ of $w$ with endpoints in two different connected components of $z$ (or vice versa)		
		\State \Call {UpdateFixedEdgesQueue}{$z,w,e,fixEdgesN$}
		\State $z',w' \gets$ \Call {GetNewCycleCovers}{$z \cup w$}	\Comment{Third neighborhood structure}
		\If {$e^{- \frac{E(z',w') - E(z,w)}{T}} \geq$ \Call{Rand}{0,1}}		\Comment{SA probabilistic transition}
		\State  $z \gets z'$, $w \gets w'$			\Comment{Proceed to a new solution}
		\EndIf
		\State $T \gets$ \Call{CoolingSchedule}{$k$}
		\State \Return $z$ and $w$
		\EndProcedure	
	\end{algorithmic}
\end{algorithm}

For the third neighborhood structure, we chose the simulated annealing approach \cite{kirk:1983} for two reasons.
Firstly, constructing a cycle cover through perfect matching is an expensive enough procedure to examine the entire neighborhood.
Secondly, the simulated annealing allows the algorithm to get out of a local minimum.

\subsection{Final algorithm}

The final general variable neighborhood search algorithm is summarized in Algorithm~\ref{Alg:GVNS}.
We have arranged the neighborhoods in order of increasing complexity:
\begin{enumerate}
	\item the first neighborhood is quite small and simple ($O(V)$ for directed graphs and $O(k \cdot V^2)$ for undirected graphs);
	
	\item the second neighborhood explores the entire search tree to a limited depth, 
	therefore, it is larger and requires significant time ($O(2^d \cdot V)$ for directed graphs and $O(3^d \cdot V)$ for undirected graphs);
	
	\item construction of each neighboring solution with respect to the third neighborhood requires a recalculation of vertex-disjoint cycle covers by calling perfect-matching algorithms, 
	therefore, it is the most time-consuming of the three. 
\end{enumerate}

We combined the first two neighborhoods in the basic variable neighborhood descent approach \cite{duar:san:mlad:tod:2018}: with the first improvement w.r.t. the second neighborhood, we return to the first neighborhood.
If the algorithm is stuck in a local minimum w.r.t. both first and second neighborhood structures, we apply the simulated annealing approach w.r.t. the third neighborhood to shake the solution and get new cycle covers.

\begin{algorithm}[t]
	\caption{General variable neighborhood search}\label{Alg:GVNS}
	\begin{algorithmic}[1]
		\Procedure{GVNS}{$x \cup y,iterN,initT,fixEdgesN,depthL,timeL$}
		\State $T \gets initT$
		\State $z,w \gets$ \Call{GetInitialCycleCovers}{$x \cup y$} \Comment{By perfect matching algorithms}
		\State $iter \gets 1$
		\While {$iter \leq iterN$ \textbf{and} running time $\leq timeL$}
		\Repeat 									\Comment{Variable neighborhood descent}
		\If {$z$ and $w$ are Hamiltonian cycles (different from $x$ and $y$)}
		\State \Return Hamiltonian decomposition $z$ and $w$
		\EndIf
		\State \Call {LocalSearch\_1} {$z,w$}	\Comment{First improvement w.r.t. the first neighborhood}
		\If {a better solution is found}		
		\State \textbf{continue}			\Comment{Restart the local search w.r.t. the first neighborhood}
		\EndIf									\Comment{Local minimum w.r.t. the first neighborhood}
		\State \Call {LocalSearch\_2} {$z,w,depthL$}	\Comment{First improvement w.r.t. the second neighborhood}
		\Until $z$ and $w$ is a local minimum w.r.t. both neighborhood structures						
		\Repeat										 \Comment{Third neighborhood structure}
		\State  \Call {SimulatedAnnealing} {$z,w,iter,T,fixEdgesN$}
		\State $iter \gets iter + 1$
		\Until {a new solution is found}
		\EndWhile
		\State \Return Hamiltonian decomposition into cycles $z$ and $w$ is not found
		\EndProcedure		
	\end{algorithmic}
\end{algorithm}

\subsection{Stopping criteria}

If the decomposition of the multigraph $x \cup y$ into edge-disjoint Hamiltonian cycles $z$ and $w$ is found, then the algorithm successfully stops and returns the solution. 
By the sufficient condition (Lemma~\ref{lemma_sufficient}), the corresponding vertices $x^v$ and $y^v$ of the traveling salesperson polytope are not adjacent.
Otherwise, we stop when the limits are exceeded either for the number of shaking steps $iterN$ or for the running time $timeL$.
The algorithm could not find the Hamiltonian decomposition, hence the corresponding vertices are ``probably adjacent''.
Thus, the algorithm has a one-sided error: the answer ``not adjacent'' is always correct, while the answer ``probably adjacent'' leaves the possibility that the vertices are not adjacent.

\subsection{Vertex adjacency and preprocessing}

The GVNS algorithm is designed primarily for the construction of the Hamiltonian decomposition of 4-regular multigraph $x \cup y$.
The only time the algorithm directly accesses the cycles $x$ and $y$ is to verify that the Hamiltonian decomposition $z$ and $w$ is different from $x$ and $y$ (line 7, Algorithm~\ref{Alg:GVNS}).
Thus, if we omit this check, then we obtain the algorithm for the Hamiltonian decomposition problem without any connection to the traveling salesperson polytope.
This is one of the main features of the algorithm.

Otherwise, if we are interested exactly in vertex adjacency in 1-skeleton of the traveling salesperson polytope, then the preprocessing step can be added to the algorithm, at which we first consider sufficient conditions for nonadjacency that can be verified in polynomial time. 
For example, pyramidal tours \cite{bond:nik:2018}, pyramidal tours with step-backs \cite{nik:2019}, pedigrees \cite{arth:2006,arth:2013}, etc.

\section{Computational results}

We chose four algorithms for comparison:
\begin{enumerate}
	\item SA: simulated annealing with respect to the third neighborhood only (this is the previous version of the algorithm described in \cite{kozl:nik:2019});
	
	\item GVNS: general variable neighborhood search algorithm presented in this paper (Algorithm~\ref{Alg:GVNS}).
	
	\item VND-12: one iteration of basic variable neighborhood descent w.r.t. the first and second neighborhood structures (lines 6-15 of the Algorithm~\ref{Alg:GVNS});
	
	\item VND-21: the same basic variable neighborhood descent algorithm as VND-12, but the order of the first and second neighborhoods is inverted.
\end{enumerate}

The considered algorithms were tested on random directed and undirected Hamiltonian cycles and on directed and undirected pyramidal tours.
A Hamiltonian tour is called \textit{a pyramidal} if the salesperson starts in city 1, then
visits some cities in ascending order, reaches city $n$, and returns to city 1 visiting
the remaining cities in descending order.
We chose the pyramidal tours since for them the Hamiltonian decomposition problem of whether for two given pyramidal tours $x$ and $y$ the union multigraph $x \cup y$ contains a pair of edge-disjoint pyramidal tours $z$ and $w$ can be solved in linear time \cite{bond:nik:2018}.
We used this fact as a sufficient condition for nonadjacency and generated random pairs of pyramidal tours for which a solution to the Hamiltonian decomposition problem is guaranteed to exist.

We used the following parameters:
\begin{itemize}
	\item for all algorithms: the time limit $timeL = 500$ seconds;
	\item for SA: the number of iterations $iterN=5\,000$, the initial temperature $initT=1\,000$, and the size of the queue of fixed edges $fixEdgesN = \lfloor \frac{n}{3}\rfloor$;
	\item for GVNS: the number of iterations $iterN=1\,000$, the recursion depth for the second neighborhood $depthL = 10$;  the initial temperature $initT=1\,000$, and the size of the queue of fixed edges $fixEdgesN = \lfloor \frac{n}{3}\rfloor$;
	\item for VND-12 and VND-21: the recursion depth for the second neighborhood $depthL = 10$.
\end{itemize}

The algorithms are implemented in Node.js. 
All experiments are performed on a machine with an Intel(R) Core(TM) i5-4690 with CPU
3.50GHz and 16GB RAM.

The results of the tests are presented in Tables~\ref{table:random_undirected}--\ref{table:pyramidal_directed}.
For each graph size, from $N=96$ to $N=6\,144$, we generated 100 random tests.
Since the Hamiltonian decomposition is a decision problem, we measured only three characteristics: the average runtime in seconds for solved and unsolved instances, and the percentage of successfully solved problems out of 100.

{\footnotesize
\begin{table}[p]
	\centering
	\caption{Computational results for 100 random undirected Hamiltonian cycles}
	\label{table:random_undirected}
	\begin{tabular}{|*{7}{r|}}
		\hline
		& \multicolumn{3}{c|}{SA} &  \multicolumn{3}{c|}{GVNS} \\ 
		\hline
		N & Not solved (s) & Solved (s) & Not adj.  & Not solved (s) & Solved (s) & Not adj.  \\
		\hline
		256 & $-$ & $2.399$ & $100\%$ & $-$ & $0.039$ & $100\%$ \\ 
		\hline
		384 & $-$ & $7.736$ & $100\%$ & $-$ & $0.077$ & $100\%$ \\
		\hline
		512 & $-$ & $13.225$ & $100\%$ & $-$ & $0.139$ & $100\%$ \\
		\hline
		768 & $-$ & $47.686$ & $100\%$ & $-$ & $0.285$ & $100\%$ \\
		\hline
		1024 & $-$ & $114.389$ & $100\%$ & $-$ & $0.564$ & $100\%$ \\
		\hline
		1536 & $-$ & $378.941$ & $100\%$ & $-$ & $1.162$ & $100\%$ \\ 
		\hline
		2048 & $500.000$ & $176.635$ & $44\%$ & $-$ & $1.901$ & $100\%$ \\
		\hline
		3072 & $500.000$ & $223.015$ & $25\%$ & $-$ & $3.896$ & $100\%$ \\
		\hline
		4096 & $500.000$ & $230.985$ & $14\%$ & $-$ & $8.265$ & $100\%$ \\
		\hline
		6144 & $500.000$ & $-$ & $0\%$ & $-$ & $17.798$ & $100\%$ \\
		\hline
		& \multicolumn{3}{c|}{VND-12} &  \multicolumn{3}{c|}{VND-21} \\ 
		\hline
		N & Not solved (s) & Solved (s) & Not adj. & Not solved (s) & Solved (s) & Not adj.  \\
		\hline
		256 & $-$ & $0.034$ & $100\%$ & $-$ & $1.919$ & $100\%$ \\
		\hline
		384 & $-$ & $0.073$ & $100\%$ & $-$ & $3.826$ & $100\%$ \\
		\hline
		512 & $-$ & $0.124$ & $100\%$ & $-$ & $4.877$ & $100\%$ \\
		\hline
		768 & $-$ & $0.254$ & $100\%$ & $-$ & $8.520$ & $100\%$ \\ 
		\hline
		1024 & $-$ & $0.459$ & $100\%$ & $-$ & $12.219$ & $100\%$ \\ 
		\hline
		1536 & $-$ & $1.151$ & $100\%$ & $-$ & $19.884$ & $100\%$ \\ 
		\hline
		2048 & $-$ & $2.055$ & $100\%$ & $-$ & $26.581$ & $100\%$ \\
		\hline
		3072 & $-$ & $5.354$ & $100\%$ & $-$ & $52.542$ & $100\%$ \\
		\hline
		4096 & $-$ & $8.660$ & $100\%$ & $-$ & $77.732$ & $100\%$ \\
		\hline
		6144 & $-$ & $18.510$ & $100\%$ & $-$ & $145.018$ & $100\%$ \\
		\hline
	\end{tabular}
\end{table}

\begin{table}[p]
	\centering
	\caption{Computational results for 100 random directed Hamiltonian cycles}
	\label{table:random_directed}
	\begin{tabular}{|*{7}{r|}}
		\hline
		& \multicolumn{3}{c|}{SA} &  \multicolumn{3}{c|}{GVNS} \\ 
		\hline
		N & Not solved (s) & Solved (s) & Not adj. & Not solved (s) & Solved (s) & Not adj.  \\
		\hline
		96 & $5.813$ & $0.125$ & $17\%$ & $9.572$ & $0.022$ & $19\%$ \\ 
		\hline
		128 & $9.143$ & $0.576$ & $17\%$ & $15.875$ & $0.130$ & $18\%$ \\  
		\hline
		192 & $16.147$ & $0.654$ & $21\%$ & $30.196$ & $0.187$ & $21\%$ \\
		\hline
		256 & $26.009$ & $1.941$ & $20\%$ & $64.056$ & $1.108$ & $25\%$ \\
		\hline
		384 & $52.432$ & $6.601$ & $17\%$ & $121.141$ & $0.933$ & $20\%$ \\ 
		\hline
		512 & $95.017$ & $19.685$ & $19\%$ & $191.424$ & $1.706$ & $22\%$ \\  
		\hline
		768 & $206.871$ & $55.055$ & $17\%$ & $395.074$ & $1.263$ & $19\%$ \\  
		\hline
		1024 & $357.549$ & $97.451$ & $14\%$ & $463.027$ & $2.153$ & $17\%$ \\ 
		\hline
		1536 & $499.232$ & $258.778$ & $9\%$ & $500.000$ & $4.972$ & $16\%$ \\ 
		\hline
		2048 & $500.000$ & $30.552$ & $6\%$ & $500.000$ & $26.161$ & $26\%$ \\ 
		\hline
		& \multicolumn{3}{c|}{VND-12} &  \multicolumn{3}{c|}{VND-21} \\ 
		\hline
		N & Not solved (s) & Solved (s) & Not adj. & Not solved (s) & Solved (s) & Not adj.  \\
		\hline
		96 & $0.098$ & $0.015$ & $7\%$ & $0.525$ & $0.472$ & $4\%$ \\
		\hline
		128 & $0.217$ & $0.003$ & $7\%$ & $1.039$ & $1.187$ & $9\%$ \\ 
		\hline
		192 & $0.293$ & $0.006$ & $10\%$ & $1.835$ & $1.257$ & $11\%$ \\ 
		\hline
		256 & $0.843$ & $0.150$ & $11\%$ & $3.172$ & $2.494$ & $9\%$ \\ 
		\hline
		384 & $1.294$ & $0.152$ & $7\%$ & $5.079$ & $5.131$ & $10\%$ \\ 
		\hline
		512 & $1.669$ & $0.151$ & $9\%$ & $7.804$ & $3.866$ & $7\%$ \\ 
		\hline
		768 & $3.079$ & $0.059$ & $4\%$ & $12.063$ & $10.049$ & $8\%$ \\ 
		\hline
		1024 & $4.888$ & $0.358$ & $6\%$ & $16.324$ & $15.265$ & $3\%$ \\ 
		\hline
		1536 & $7.772$ & $0.177$ & $5\%$ & $31.842$ & $33.392$ & $5\%$ \\ 
		\hline
		2048 & $9.602$ & $0.260$ & $6\%$ & $57.515$ & $48.423$ & $4\%$ \\ 
		\hline
	\end{tabular}
\end{table}

\begin{table}[p]
	\centering
	\caption{Computational results for 100 undirected pyramidal tours}
	\label{table:pyramidal_undirected}
	\begin{tabular}{|*{7}{r|}}
		\hline
		& \multicolumn{3}{c|}{SA} &  \multicolumn{3}{c|}{GVNS} \\ 
		\hline
		N & Not solved (s) & Solved (s) & Not adj. & Not solved (s) & Solved (s) & Not adj.  \\
		\hline
		256 & $29.291$ & $7.362$ & $7\%$ & $-$ & $0.062$ & $100\%$ \\ 
		\hline
		384 & $60.223$ & $18.904$ & $2\%$ & $-$ & $0.165$ & $100\%$ \\
		\hline
		512 & $106.710$ & $-$ & $0\%$ & $-$ & $0.252$ & $100\%$ \\
		\hline
		768 & $235.904$ & $-$ & $0\%$ & $-$ & $0.623$ & $100\%$ \\ 
		\hline
		1024 & $376.132$ & $-$ & $0\%$ & $-$ & $1.129$ & $100\%$ \\ 
		\hline
		1536 & $500.000$ & $-$ & $0\%$ & $-$ & $2.319$ & $100\%$ \\ 
		\hline
		2048 & $500.000$ & $-$ & $0\%$ & $-$ & $4.322$ & $100\%$ \\
		\hline
		3072 & $500.000$ & $-$ & $0\%$ & $-$ & $12.933$ & $100\%$ \\ 
		\hline
		4096 & $500.000$ & $-$ & $0\%$ & $-$ & $30.939$ & $100\%$ \\
		\hline
		6144 & $500.000$ & $-$ & $0\%$ & $-$ & $52.279$ & $100\%$ \\
		\hline
		& \multicolumn{3}{c|}{VND-12} &  \multicolumn{3}{c|}{VND-21} \\ 
		\hline
		N & Not solved (s) & Solved (s) & Not adj. & Not solved (s) & Solved (s) & Not adj. \\
		\hline
		256 & $0.155$ & $0.074$ & $97\%$ & $0.291$ & $0.544$ & $96\%$ \\ 
		\hline
		384 & $-$ & $0.155$ & $100\%$ & $1.531$ & $0.882$ & $99\%$ \\ 
		\hline
		512 & $0.473$ & $0.253$ & $97\%$ & $2.943$ & $1.731$ & $96\%$ \\
		\hline
		768 & $2.212$ & $0.609$ & $96\%$ & $5.529$ & $3.732$ & $98\%$ \\ 
		\hline
		1024 & $-$ & $1.059$ & $100\%$ & $6.040$ & $7.468$ & $96\%$ \\ 
		\hline
		1536 & $8.723$ & $3.004$ & $98\%$ & $15.815$ & $16.650$ & $95\%$ \\ 
		\hline
		2048 & $-$ & $5.730$ & $100\%$ & $46.727$ & $34.806$ & $97\%$ \\ 
		\hline
		3072 & $-$ & $13.488$ & $100\%$ & $94.374$ & $87.131$ & $92\%$ \\ 
		\hline
		4096 & $73.669$ & $25.677$ & $99\%$ & $280.352$ & $201.532$ & $94\%$ \\ 
		\hline
		6144 & $-$ & $60.268$ & $100\%$ & $500.000$ & $285.637$ & $32\%$ \\
		\hline
	\end{tabular}
\end{table}

\begin{table}[p]
	\centering
	\caption{Computational results for 100 directed pyramidal tours}
	\label{table:pyramidal_directed}
	\begin{tabular}{|*{7}{r|}}
		\hline
		& \multicolumn{3}{c|}{SA} &  \multicolumn{3}{c|}{GVNS} \\ 
		\hline
		N & Not solved (s) & Solved (s) & Not adj. & Not solved (s) & Solved (s) & Not adj.  \\
		\hline
		192 & $12.066$ & $2.462$ & $72\%$ & $-$ & $0.034$ & $100\%$ \\ 
		\hline
		256 & $18.492$ & $3.324$ & $45\%$ & $-$ & $0.078$ & $100\%$ \\ 
		\hline
		384 & $35.991$ & $4.426$ & $24\%$ & $-$ & $0.224$ & $100\%$ \\ 
		\hline
		512 & $59.528$ & $7.112$ & $15\%$ & $-$ & $0.582$ & $100\%$ \\ 
		\hline
		768 & $120.154$ & $13.237$ & $8\%$ & $-$ & $1.820$ & $100\%$ \\  
		\hline
		1024 & $202.288$ & $-$ & $0\%$ & $-$ & $3.001$ & $100\%$ \\
		\hline
		1536 & $269.308$ & $-$ & $0\%$ & $500.000$ & $11.000$ & $96\%$ \\ 
		\hline
		2048 & $374.129$ & $-$ & $0\%$ & $500.000$ & $22.805$ & $95\%$ \\
		\hline
		3072 & $500.000$ & $-$ & $0\%$ & $500.000$ & $32.614$ & $83\%$ \\
		\hline
		4096 & $500.000$ & $-$ & $0\%$ & $500.000$ & $72.210$ & $70\%$ \\
		\hline
		& \multicolumn{3}{c|}{VND-12} &  \multicolumn{3}{c|}{VND-21} \\ 
		\hline
		N & Not solved (s) & Solved (s) & Not adj. & Not solved (s) & Solved (s) & Not adj.  \\
		\hline
		192 & $0.020$ & $0.009$ & $73\%$ & $0.065$ & $0.105$ & $66\%$ \\ 
		\hline
		256 & $0.043$ & $0.017$ & $74\%$ & $0.126$ & $0.170$ & $67\%$ \\ 
		\hline
		384 & $0.106$ & $0.028$ & $66\%$ & $0.385$ & $0.421$ & $53\%$ \\
		\hline
		512 & $0.152$ & $0.064$ & $56\%$ & $0.645$ & $0.870$ & $61\%$ \\ 
		\hline
		768 & $0.322$ & $0.132$ & $61\%$ & $1.827$ & $1.914$ & $44\%$ \\ 
		\hline
		1024 & $0.580$ & $0.255$ & $46\%$ & $3.428$ & $3.311$ & $36\%$ \\
		\hline
		1536 & $1.209$ & $0.470$ & $38\%$ & $9.865$ & $9.433$ & $30\%$ \\ 
		\hline
		2048 & $2.058$ & $0.957$ & $36\%$ & $22.995$ & $20.554$ & $27\%$ \\
		\hline
		3072 & $5.483$ & $3.087$ & $17\%$ & $66.456$ & $53.724$ & $14\%$ \\
		\hline
		4096 & $10.399$ & $3.984$ & $16\%$ & $193.328$ & $143.103$ & $12\%$ \\
		\hline
	\end{tabular}
\end{table}
}

Based on the test results, it can be seen that for random undirected Hamiltonian cycles (Table~\ref{table:random_undirected}) all algorithms proved to be quite successful.
Although, the GVNS on average turned out to be 150 times faster than SA, which allowed GVNS to solve all instances, while SA on graphs with more than 2\,000 vertices ran into a time limit and failed.
Note that a single iteration of VND was enough to solve all instances, with VND-12 working on average 25 times faster than VND-21.
We can conclude that the union of two random undirected Hamiltonian cycles almost always contains several different Hamiltonian decompositions and any of the algorithms finds at least one of them.
From a geometric point of view, two random vertices of the symmetric traveling salesperson polytope $\TSP(n)$ are not adjacent with a very high probability.

The situation is completely different for random directed Hamiltonian cycles (Table~\ref{table:random_directed}),
where the SA and GVNS found a solution in about 10-25\% of cases.
We investigated the behavior of algorithms on small graphs (up to 50 vertices), where the problem can be solved by backtracking, and the SA and GVNS showed accurate results there.
From a geometric perspective, it is known that the $1$-skeleton of the asymmetric traveling salesperson polytope is generally much denser than the $1$-skeleton of the symmetric polytope. For example, the diameter of the $1$-skeleton of $\ATSP(n)$ is $2$ \cite{padb:rao:1974}, while the best known upper bound for the diameter of the $1$-skeleton of $\TSP(n)$ is $4$ \cite{risp:1998}.

Summary for random directed Hamiltonian cycles: VND-21 solved 70 instances, VND-12 solved 72 instances, SA solved 157 instances and GVNS solved 203 instances out of $1\,000$.
It can be seen that for directed graphs VND w.r.t. the first and second neighborhood structures is not as effective as for undirected graphs.
Between two versions of variable neighborhood descent, the VND-12 was on average 5 times faster.
We compared the SA and GVNS by their success rate.
Note that it is hidden in the Table~\ref{table:random_directed} that their results are nested: all instances solved by SA were also solved by GVNS.
Since we are considering the decision problem with yes/no answers, we applied the McNemar test to evaluate the statistical significance of the results \cite{mcnem:1947}. 
The GVNS algorithm against SA had $\chi^2$ with Yates correction of 1.0 equal to $44.022$ with the corresponding $p$-value being less than $1 \times 10^{-6}$.
Thus, the GVNS algorithm showed statistically significantly better results for random directed Hamiltonian cycles.

However, in the general case, since we tested the algorithms on random cycles, we do not know if the algorithm could not find a solution because it failed or because there is no solution.
Therefore, we considered the synthetic tests on pyramidal tours, for which the sufficient condition for nonadjacency is satisfied and the solution is guaranteed to exist.

Test results on pyramidal tours (Tables \ref{table:pyramidal_undirected}-\ref{table:pyramidal_directed}) allow even better comparison between the GVNS and the previous SA version of the algorithm.
The SA performed poorly and solved only 9 instances for undirected graphs and 164 instances for directed graphs out of $1\,000$, while the GVNS solved all $1\,000$ instances for undirected graphs and 944 instances for directed graphs.
Besides, the GVNS was much faster. 
Again, it can be seen that the variable neighborhood descent is very successful for undirected graphs where just one iteration solved almost all instances (987 for VND-12 and 895 for VND-21), and slightly less effective for directed graphs (483 instances for VND-12 and 410 for VND-21 out of $1\,000$).
Between two versions of VND, their success rate was relatively close with VND-12 being on average 6 times faster for undirected graphs and 10 times faster for directed graphs.

Summing up the tests, the general variable neighborhood search showed significant improvement both in success rate and running time over the simulated annealing approach from \cite{kozl:nik:2019}. 
For undirected graphs, often just one iteration of the variable neighborhood descent w.r.t. the first and second neighborhood structures was enough.
The only weak point of GVNS was on directed pyramidal tours of large size, where the algorithm failed to solve several instances with a solution known to exist. However, each time the algorithm ran into a time limit. With broader time constraints, the results could be better.
Also note that the alternative order of the neighborhoods, different from the increasing complexity, only slowed down the algorithm and did not lead to better results.

\section{Conclusion} 

We introduce the general variable neighborhood search algorithm with 3 different neighborhood structures to find a Hamiltonian decomposition of the 4-regular multigraph. The algorithm showed good computational results on both directed and undirected graphs.

The Hamiltonian decomposition problem arises in the analysis of the 1-skeleton of the traveling salesperson polytope.
Thus, the GVNS algorithm can be applied to verify vertex adjacency in the 1-skeleton of both $\TSP(n)$ and $\ATSP(n)$ polytopes based on the sufficient condition of Lemma~\ref{lemma_sufficient}.
Here the algorithm has a one-sided error: the answer ``not adjacent'' is always correct, while the answer ``probably adjacent'' leaves the possibility that the vertices actually are not adjacent.

With some modifications, the algorithm can be applied to the Hamiltonian decomposition problem in $2k$-regular graphs.

\subsection*{Acknowledgements}
We are very grateful to the anonymous reviewers for their comments and suggestions which helped to improve the
presentation of the results in this paper.

\end{document}